\documentclass[12pt]{article}
\usepackage{amsfonts}
\def\mod{\mathop{\mathrm{mod}}}
\def\Box{\raisebox{1pt}{\framebox{\phantom{\small x}}}}
\newtheorem{thm}{Theorem}
\newtheorem{cor}{Corollary}
\newenvironment{prf}{\begin{trivlist}\item[]{\bf Proof} }%
{\hfill $\Box$ \end{trivlist}}
\begin{document}
\renewcommand{\thefootnote}{$\sharp$}
\renewcommand{\section}{\subsection}
\begin{center}
{\bf\Large Some Special Geometry in Dimension Six}\\[10pt]
{\bf\large Andreas \v Cap and Michael Eastwood\footnote{Senior Research Fellow
of the Australian Research Council.}}
\end{center}

\noindent {\small {\bf Abstract:}\quad We generalise the notion of contact
manifold by allowing the contact distribution to have codimension two.
There are special features in dimension six. In particular, we show that the
complex structure on a three-dimensional complex contact manifold is determined
solely by the underlying contact distribution.}

\renewcommand{\thefootnote}{}
\footnotetext{This research was undertaken whilst the second author was
visiting the Erwin Schr\"odinger International Institute
for Mathematical Physics. Its support is gratefully acknowledged.}

\section*{Definitions}
Suppose $M$ is a 6-dimensional connected oriented
smooth manifold and $H$
is a rank 4 smooth subbundle of its
tangent bundle~$TM$. Let $Q$ denote the quotient bundle~$TM/H$. There is
a homomorphism of vector bundles ${\mathcal L}:H\wedge H\to Q$
induced by Lie bracket:--
$${\mathcal L}(\xi,\eta)=[\xi,\eta]\mod H\quad
                                         \mbox{for }\xi,\eta\in\Gamma(H).$$
Regard ${\mathcal L}$ as a tensor
${\mathcal L}\in\Gamma(\Lambda^2H^*\otimes Q)$. Then
${\mathcal L}\wedge{\mathcal L}\in \Gamma(\Lambda^4H^*\otimes\bigodot^2Q)$
may be regarded as a quadratic form on $Q^*$ defined up to scale.
We shall say that $(M,H)$ is {\em non-degenerate\/}
if and only if
${\mathcal L}\wedge{\mathcal L}$ is non-degenerate as such a quadratic form.
Since $Q$ has rank two, there are only two cases:--
\begin{itemize}
\item $(M,H)$ is {\em elliptic\/}
      $\iff {\mathcal L}\wedge{\mathcal L}$ is definite;
\item $(M,H)$ is {\em hyperbolic\/}
      $\iff {\mathcal L}\wedge{\mathcal L}$ is indefinite.
\end{itemize}
An elliptic example may be obtained by taking a
3-dimensional complex contact manifold and forgetting its complex structure.
A hyperbolic example may be obtained by taking the product
of two 3-dimensional real contact manifolds. These two examples will be
referred to as the `flat' models. The motivations for our investigation are
discussed in the end of this article.

\section*{Acknowledgements}
We are pleased to
acknowledge several useful conversations with Gerd Schmalz, Jan Slov\'ak, and
Peter Vassiliou.

\section*{The Elliptic Case}
\begin{thm} Suppose $(M,H)$ is elliptic. Then $M$ admits a unique almost
complex structure $J:TM\to TM$ characterised by the following properties:--
\begin{itemize}
\item $J$ preserves $H$;
\item the orientation on $M$ induced by $J$ is the given one;
\item $\mathcal L:H\times H\to Q$ is complex bilinear for the induced
structures, or equivalently $[\xi,\eta]+J[J\xi,\eta]\in\Gamma(H)\quad
              \mbox{for }\xi,\eta\in\Gamma(H)$;
\item $[\xi,\eta]+J[J\xi,\eta]-J[\xi,J\eta]+[J\xi,J\eta]\in\Gamma(H)\quad
              \mbox{for }\xi\in\Gamma(TM),\;\eta\in\Gamma(H).$
\end{itemize}
Furthermore, the tensor $S:Q\otimes H\to Q$ induced by
$$S(\xi,\eta)=[\xi,\eta]+J[J\xi,\eta]\mod H\quad
              \mbox{for }\xi\in\Gamma(TM),\;\eta\in\Gamma(H)$$
is the obstruction to $J$ being integrable.
\end{thm}
\begin{prf}
Fix $x\in M$.
Since ${\mathcal L}_x\wedge{\mathcal L}_x$ is definite, there is no
$\psi\in Q_x^*$ for which
$(\psi\circ{\mathcal L}_x)\wedge(\psi\circ{\mathcal L}_x)$ vanishes---as a
quadratic polynomial ${\mathcal L}_x\wedge{\mathcal L}_x$ has no real roots.
Instead it has two complex roots, related by complex conjugation.
Each of these roots gives $\psi\in Q_x^*\otimes{\mathbb C}$ defined up to
complex scale, so that
$(\psi\circ{\mathcal L}_x)\wedge(\psi\circ{\mathcal L}_x)$
vanishes as an element of $\Lambda^4H^*\otimes{\mathbb C}$. In this case,
according to the
Pl\"ucker criterion, $\psi\circ{\mathcal L}_x$ is simple as an element of
$\Lambda^2H^*\otimes{\mathbb C}$. The corresponding complex 2-plane in
$H^*\otimes{\mathbb C}$ defines a complex structure $J:H_x\to H_x$.
At the same time $\psi\in Q_x^*\otimes{\mathbb C}$ identifies $Q_x$ with
${\mathbb C}$ and, in particular, defines a complex structure $J:Q_x\to Q_x$.
These complex structures on $H_x$ and $Q_x$ are unchanged if $\psi$ is
multiplied by any complex number. In other words, they are determined by
choosing one of the two roots of ${\mathcal L}_x\wedge{\mathcal L}_x$
as a quadratic polynomial. The other root replaces $J$ by $-J$ but only one of
these choices induces the given orientation on~$M$. To summarise, we now have
uniquely determined almost complex structures on $H$ and $Q$ so that
\begin{equation}\label{JL}{\mathcal L}(\xi,\eta)+J{\mathcal L}(J\xi,\eta)=0
  \quad\mbox{for }\xi,\eta\in\Gamma(H)\end{equation}
and inducing the given orientation on~$M$.
Choose any extension of these almost complex structures to an almost
complex structure $\tilde J:TM\to TM$. This $\tilde J$ satisfies the
first three properties claimed in the
statement of the theorem. Define $\tilde S:TM \otimes H\to Q$ by
\begin{equation}\label{defS}\tilde
S(\xi,\eta)=[\xi,\eta]+J[\tilde J\xi,\eta]\mod H\quad
\mbox{for }\xi\in\Gamma(TM),\;\eta\in\Gamma(H)\end{equation}
This homomorphism depends on the choice of the etension $\tilde
J$. For fixed $\xi\in TM$ consider the map $H\to Q$ defined by
$\eta\mapsto\frac{1}{2}(-\tilde S(\xi,\eta)+J\tilde S(\xi,J\eta))$. By
construction this map is complex linear, so
non--degeneracy of $\mathcal L$ implies that there is a unique element
$K\xi\in H$ such that
\begin{equation}\label{defK}
{\mathcal L}(K\xi,\eta)=\frac{-\tilde S(\xi,\eta)+J\tilde
S(\xi,J\eta)}{2}\quad \mbox{for
}\xi\in\Gamma(TM),\;\eta\in\Gamma(H)\end{equation}
and this defines a homomorphism $K:TM\to H$.

We claim that $J=\tilde J+K$ is the almost complex structure whose
existence is asserted in the statement of the theorem.
If $\xi\in\Gamma(H)$, then (\ref{JL}) implies that
$\tilde S(\xi,\eta)=0$ so $K\xi=0$, and in particular
$K^2=0$. Therefore, $J$ preserves~$H$. Also
$$(\tilde J+K)^2=\tilde J^2+\tilde JK+K\tilde J+K^2=-\mbox{Id}+\tilde
JK+K\tilde J$$
so we must check that $\tilde JK+K\tilde J=0$. By the non-degeneracy
of ${\mathcal L}$ it suffices to check that
$${\mathcal L}(\tilde JK\xi,\eta)+{\mathcal L}(K\tilde J\xi,\eta)=0\quad
\mbox{for }\xi\in\Gamma(TM),\;\eta\in\Gamma(H).$$
This is easily verified using (\ref{JL}), (\ref{defS}), and~(\ref{defK}).
Thus, $J$ is an almost complex structure. Moreover, it
satisfies the first three
requirements listed in the theorem as a consequence of $\tilde J$
doing so. Moreover, the tensor $S$ corresponding to $J=\tilde J+K$ is
visibly given by $S(\xi,\eta)=\tilde S(\xi,\eta)+\mathcal
L(K\xi,\eta)$. By construction, this is just the component of $\tilde S$
which is conjugate linear in the second variable. But the final
requirement is immediately seen to be equivalent to the fact that the
corresponding tensor $S$ (which is conjugate linear in the first
variable by construction), is conjugate linear in the second variable,
too. In fact, this forces (\ref{defK}) as the correct modification so
$J$ is uniquely characterised by having all four properties.

It remains to show that the tensor $S$ is the obstruction to
integrability of $J$. The Nijenhuis tensor of $J$ is
$$N(\xi,\eta)=[\xi,\eta]+J[J\xi,\eta]+J[\xi,J\eta]-[J\xi,J\eta]\quad
\mbox{for }\xi,\eta\in\Gamma(TM).$$
Notice that $N$ is skew and $N(\xi,J\eta)=-JN(\xi,\eta)$. In particular,
\begin{equation}\label{NJ}N(\xi,J\xi)=-JN(\xi,\xi)=0\quad
\mbox{for }\xi\in\Gamma(TM).\end{equation}
Firstly, consider the case when $\xi\in\Gamma(TM),\;\eta\in\Gamma(H)$.
The vanishing of $S$ means that
\begin{equation}\label{ass}[\xi,\eta]+J[J\xi,\eta]\in\Gamma(H)\quad
\mbox{for }\xi\in\Gamma(TM),\;\eta\in\Gamma(H).\end{equation}
In particular, this implies $N(\xi,\eta)\in\Gamma(H)$, so
we may consider the tensor $R:TM\otimes H\otimes H\to Q$ defined by
$$R(\xi,\eta,\mu)={\mathcal L}(N(\xi,\eta),\mu)\quad
\mbox{for }\xi\in\Gamma(TM),\;\eta,\mu\in\Gamma(H).$$
We claim that $R$ vanishes. Once this is proved, non-degeneracy of
${\mathcal L}$ implies that
$N(\xi,\eta)=0$ for $\xi\in\Gamma(TM),\;\eta\in\Gamma(H)$ and so $N$ descends
to $N:\Lambda^2Q\to TM$.
Then, as $Q$ has complex rank one, (\ref{NJ}) forces $N$ to vanish.

To complete the proof, therefore, it suffices to show that $R$ vanishes. In the
following calculation $\equiv$ denotes equality modulo~$H$ and in passing from
one line to the next we are using either the Jacobi identity,
or~(\ref{ass}), or the fact that $S$ is conjugate linear in both
variables.
$$\begin{array}{rcl}R(\xi,\eta,\mu)
&\equiv&[[\xi,\eta],\mu]+[J[J\xi,\eta],\mu]
+[J[\xi,J\eta],\mu]-[[J\xi,J\eta],\mu]\\[5pt]
&\equiv&[[\xi,\eta],\mu]+J[[J\xi,\eta],\mu]
+J[[\xi,J\eta],\mu]-[[J\xi,J\eta],\mu]\\[5pt]
&=&[[\xi,\mu],\eta]+J[[J\xi,\mu],\eta]
+J[[\xi,\mu],J\eta]-[[J\xi,\mu],J\eta]\\[2pt]
&&+[[\mu,\eta],\xi]+J[[\mu,\eta],J\xi]
+J[[\mu,J\eta],\xi]-[[\mu,J\eta],J\xi]\\[5pt]
&\equiv&[[\xi,\mu],\eta]+[J[J\xi,\mu],\eta]
+[J[\xi,\mu],J\eta]-[[J\xi,\mu],J\eta]\\[2pt]
&&+[[\mu,\eta],\xi]+J[[\mu,\eta],J\xi]
+J[[\mu,J\eta],\xi]-[[\mu,J\eta],J\xi]\\[5pt]
&=&[[\xi,\mu]+J[J\xi,\mu],\eta]
+[J[\xi,\mu]-[J\xi,\mu],J\eta]\\[2pt]
&&+[[\mu,\eta],\xi]+J[[\mu,\eta],J\xi]
+J[[\mu,J\eta],\xi]-[[\mu,J\eta],J\xi]\\[5pt]
&\equiv&[[\mu,\eta],\xi]+J[[\mu,\eta],J\xi]
+J[[\mu,J\eta],\xi]-[[\mu,J\eta],J\xi].
\end{array}$$
Therefore,
$$R(\xi,\eta,\mu)+R(\xi,\mu,\eta)\equiv
J[[\mu,J\eta]+[\eta,J\mu],\xi]-[[\mu,J\eta]+[\eta,J\mu],J\xi]$$
and since $[\mu,J\eta]+[\eta,J\mu]\in\Gamma(H)$, this expression
vanishes by~(\ref{ass}).
We conclude that $R:TM\otimes H\otimes H\to Q$ is skew in its last two
entries. But by definition $R$ is conjugate linear in the middle
variable and complex linear in the last variable, which together with
skew symmetry in these two variables forces $R$ to vanish as required.
\end{prf}
\begin{cor}
The only local invariant of an elliptic $(M,H)$ is the tensor~$S$.
\end{cor}
\begin{prf} If $S$ vanishes, then $(M,H)$ is a complex contact
manifold. The Darboux theorem in the holomorphic setting says that all
3-dimensional complex contact manifolds are locally isomorphic.
\end{prf}

\section*{The Hyperbolic Case}

There is an entirely parallel story for the hyperbolic case with almost complex
structure replaced by almost product structure. The corresponding theorem may
be stated as follows.
\begin{thm} Suppose $(M,H)$ is hyperbolic. Then $H$ admits a canonical splitting
$H=H_+\oplus H_-$ characterised by the following properties:--
\begin{itemize}
\item $[\xi,\eta]\in\Gamma(H)\quad\mbox{for }
              \xi\in\Gamma(H_+),\;\eta\in\Gamma(H_-)$;
\item the orientation on M induced by
      $\xi_1\wedge\xi_2\wedge[\xi_1,\xi_2]\wedge
       \eta_1\wedge\eta_2\wedge[\eta_1,\eta_2]$ for
      $\xi_1,\xi_2\in\Gamma(H_+),\;\eta_1,\eta_2\in\Gamma(H_-)$
      is the given one.
\end{itemize}
Let $Q_\pm$ be the range of ${\mathcal L}|_{\Lambda^2H_\pm}$. Non-degeneracy of
${\mathcal L}$ implies that $Q=Q_+\oplus Q_-$.
By setting
$T_\pm M=[H_\pm,H_\pm]$,
we obtain a canonical splitting
$TM=T_+M\oplus T_-M$ such that $Q_\pm=T_\pm M/H_\pm$. Furthermore, the tensors
$S_+:Q_+\otimes H_+\to Q_-$ and $S_-:Q_-\otimes H_-\to Q_+$ induced by
$$S_\pm(\xi,\eta)\equiv[\xi,\eta]\mod (T_\pm\oplus H_\mp)\quad
  \mbox{for }\xi\in\Gamma(T_\pm M)\;,\eta\in\Gamma(H_\pm)$$
are the respective obstructions to $T_+M$ and $T_-M$ being Frobenius
integrable.
\end{thm}
If $S_\pm$ both vanish, then locally we obtain the flat local model, namely
a product of two 3-dimensional real contact
manifolds. The Darboux theorem, applied to each such contact manifold
separately, implies that the flat model is locally unique. Again, the tensors
$S_\pm$ provide the only local structure.

\section*{Motivations}

Our motivation for this article comes from the theory of CR~submanifolds of
codimension 2 in~${\mathbb C}^4$.  This theory was pioneered by Loboda~\cite{l}
and Ezhov-Schmalz~\cite{es} who found normal forms for such submanifolds
paralleling the Moser normal form for CR~hypersurfaces. In this
context, the distribution $H$ is formed by the maximal complex
subspaces of the tangent spaces. More generally, to make an elliptic
or hyperbolic $(M,H)$ into a partially integrable almost CR manifold,
one has to specify an almost complex structure $\tilde J$ on $H$ such
that $\mathcal L(\tilde J\xi,\tilde J\eta)=\mathcal L(\xi,\eta)$ for all
$\xi,\eta\in H$. In the hyperbolic case, this implies in particular
that $H=H_+\oplus H_-$ is a decomposition of $H$ as a sum of two
complex line bundles. On the other hand, in the elliptic case the
second almost complex structure $\tilde J$ can also be rephrased as a
decomposition $H=H_+\oplus H_-$ as a sum of complex line bundles
characterized by $\tilde J=\pm J$ on $H_\pm$.

Clearly, these additional structures lead to additional obstructions
against being CR--isomorphic to the flat models (which are just
appropriate quadrics). For example, one has the Nijenhuis tensor
corresponding to $\tilde J$, or the obvious obstructions against
integrability of the subbundles $H_\pm$ in the elliptic case. But in
fact, in the CR setting, one gets much more structure: In \cite{ss} and
\cite{cs} it is shown that one gets a parabolic geometry parallel to
the Chern-Moser-Tanaka theory for CR~hypersurfaces and thus in
particular canonical Cartan connections. This article may be viewed as
some remnant of the parabolic theory.

As pointed out to us by Peter Vassiliou, there is another context in
which $(M,H)$ with these special dimensions arise. The general pair of smooth
first order partial differential equations in two independent variables $(x,y)$
and two dependent variables $(u,v)$ may be regarded as a codimension 2
submanifold $M$ in the 8-dimensional jet space with co\"ordinates
$(x,y,u,v,u_x,u_y,v_x,v_y)$. This jet space has a natural distribution of rank
6 defined as common kernel of the two 1-forms
$$du-u_x\,dx-u_y\,dy\quad\mbox{and}\quad dv-v_x\,dx-v_y\,dy.$$
Generically, $M$ will meet this distribution transversally and so will
itself inherit a rank 4 distribution~$H$.
The elliptic flat model is obtained from the Cauchy-Riemann equations
$$u_x=v_y\quad\mbox{and}\quad u_y=-v_x.$$
The hyperbolic flat model arises from the equations
$$u_y=0\quad\mbox{and}\quad v_x=0.$$
Further discussion may be found in
\cite[Chapter~VII,~\S 1]{bcggg}, \cite{V}, and~\cite{v}.

\section*{Higher Dimensions}If we start with a
$(2n+1)$-dimensional complex contact manifold $M$ with contact
distribution~$H$, then ${\mathcal L}\in\Gamma(\Lambda^2H^*\otimes Q)$ may be
defined as before but now we should consider
${\mathcal L}^{\wedge 2n}\in\Gamma(\Lambda^{4n}H^*\otimes\bigodot^{2n}Q)$ as a
polynomial of degree $2n$ defined up to scale. Only when $n=1$ is this
polynomial generic. In general it has only two roots, each complex and of
multiplicity~$n$.

\begin{tabular}[t]l
{\sc Andreas \v Cap}\\
Institut f\"ur Mathematik, Universit\"at Wien\\
Strudlhofgasse~4, A-1090 Wien, Austria\\
and\\
Erwin Schr\"odinger International Institute for Mathematical Physics,\\
Boltzmanngasse~9, A-1090 Wien, Austria\\[5pt]
E-mail: andreas.cap@esi.ac.at\\[10pt]
{\sc Michael Eastwood}\\
Department of Pure Mathematics\\
University of Adelaide\\
South Australia 5005\\[5pt]
E-mail: meastwoo@maths.adelaide.edu.au
\end{tabular}
\end{document}